\documentclass[12pt]{article}
\usepackage{bbm}
\usepackage{amsfonts}
\usepackage{pifont}
\usepackage{hyperref}
\usepackage{mathrsfs}
\usepackage{indentfirst}
\usepackage{amsmath}
\usepackage{amssymb}
\usepackage[amsmath, thmmarks]{ntheorem}
\usepackage{cite}
\usepackage{graphicx}
\usepackage{tikz}
\newtheorem{definition}{\bf Definition}[section]
\newtheorem{lemma}{\bf Lemma}[section]
\newtheorem{theorem}{\bf Theorem}[section]
\newtheorem{remark}{\bf Remark}[section]

\newtheorem{proposition}{\bf Proposition}[section]
\newtheorem{example}{\bf Example}[section]

\def\QEDopen{{\setlength{\fboxsep}{0pt}\setlength{\fboxrule}{0.2pt}\fbox{\rule[0pt]{0pt}{1.3ex}\rule[0pt]{1.3ex}{0pt}}}} 
\def\QED{\QEDopen}
\def\proof{{\bf Proof.} }
\def\endproof{\hspace*{\fill}~\QED\par\endtrivlist\unskip}

\begin{document}
\setcounter{page}{1}

\title{{\textbf{Left-continuous pseudo-t-norms on modular lattices}}\thanks {Supported by the National Natural Science
Foundation of China (No. 12471440)}}
\author{Peng He\footnote{\emph{E-mail address}: 443966297@qq.com}, Xue-ping Wang\footnote{Corresponding author. xpwang1@hotmail.com; fax: +86-28-84761502}\\
\emph{\small (School of Mathematical Sciences, Sichuan Normal University}\\ \emph{\small Chengdu 610066, Sichuan, People's Republic of China)}}
\newcommand{\pp}[2]{\frac{\partial #1}{\partial #2}}
\date{}
\maketitle

\begin{quote}
{\bf Abstract}
This article focuses on the relationship between pseudo-t-norms and the structure of lattices. First, we establish a necessary and sufficient condition for the existence of a left-continuous t-norm on the ordinal sum of two disjoint complete lattices. Then, we define the $1$-distributivity of a lattice, which is applied for characterizing a complete atomistic lattice that has a left-continuous pseudo-t-norm. We also describe the forbidden structures of a finite modular lattice that is a $1$-distributive lattice, which is used for representing a kind of finite planar modular lattices that have left-continuous pseudo-t-norms.


\end{quote}

\section{Introduction}\label{intro}
Since 1970s, Zadeh \cite{Zad} and other scholars have successively introduced a series of logical operators on the 
interval $[0,1]$ to represent the relationship of `and', `or' and `not'. Bellman and Giertz \cite{Bellma}, and Alsina, Trillas and Valverde \cite{Alsina}
have conducted in-depth research on these logical operators. In fuzzy logic, connectives `and', `or' and `not' on $[0, 1]$ are typically modelled 
by t-norms, t-conorms and strong negations, respectively \cite{Weber}. Based on these logical operators on $[0,1]$, three fundamental classes of fuzzy implications, 
namely, S-implication, R-implication and QL-implication have been defined and thoroughly studied \cite{Du, Ma, Wang02, Wang03}. 

It is well known that a lattice is a natural extension of $[0,1]$ and there is a close link between 
fuzzy set theory and order theory \cite{De}. Several researchers have explored 
t-norms on bounded lattices \cite{De, Sam08, Sam, Ma, Coomann94}, for example, Ma and Wu \cite{Ma}
examined the relations between t-norms and implications on a complete lattice, 
revealing a one-to-one correspondence between the set of all 
left-continuous t-norms and the set of all right-continuous implications. 
Since properties of an implication depend on its relations with other logical operators, they also presented a 
criterion that helps us infer whether a group of logical operators is consistent. Fodor \cite{Fod} introduced a notion of a weak t-norm on $[0,1]$ 
and studied the relations between weak t-norms and fuzzy implications. Wang and Yu \cite{Wang02} investigated pseudo-t-norms and 
implication operators on a complete Brouwerian lattice $L$, showing the relations between all left-continuous pseudo-t-norms on $L$ and all right-continuous implications. 
In particular, Wang and Yu \cite{Wang03} further studied the pseudo-t-norms and implication operators on a complete Brouwerian lattice, 
and discussed their direct products and direct product decompositions. 

It should be explicitly noted that all results presented in \cite{Ma, Wang02, Wang03} are based on the existence 
of left-continuous t-norms (resp. pseudo-t-norms) on a lattice when exploring the relationship between logical operators and implication operators. 
Although left-continuous t-norms (resp. pseudo-t-norms) are guaranteed to exist on a complete Brouwerian lattice, for an arbitrarily given complete lattice, the question of whether a left-continuous t-norm exists on it remains open. 
Consequently, the problem regarding the existence of left-continuous t-norms (resp. pseudo-t-norms) on a complete lattice is in urgent need of resolution.
Additionally, what are the specific structural characterizations of the lattice in which left-continuous t-norms (resp. pseudo-t-norms) exist? This article will attempt to solve the above two problems.

The rest of this article is organized as follows. In Section 2, we recall some basic concepts and results required lately. In Section 3, we supply a necessary and sufficient condition for the existence of a left-continuous t-norm on the ordinal sum of two disjoint complete lattices. In Section 4, we define the $1$-distributivity of a lattice and utilize it to  characterize a complete atomistic lattice on which there is a left-continuous pseudo-t-norm. We also describe the forbidden structures of a finite modular lattice that is a $1$-distributive lattice. In Section 5, we equivalently describes a kind of planar $1$-distributive modular lattices on which there is a left-continuous pseudo-t-norm. A concluding remark is drawn in Section 6.

\section{Preliminaries}
For the purpose of reference, we shall briefly review some concepts and facts about lattices \cite{Birkhoff73, Crawley73, Gratz}.

A poset is a system consisting of a non-empty set $P$ and a binary relation $\leq$ in $P$ such that the following conditions are satisfied for all $x, y, z\in P$:

(1) $x\leq x$; (Reflexive)

(2) if $x\leq y$ and $y\leq x$, then $x=y$; (Antisymmetry)

(3) if $x\leq y$ and $y\leq z$, then $x\leq z$. (Transitivity)\\
The relation $\leq$ is a partial order in the set $P$. Although the set $P$ is only the domain of the partial order $\leq$, we
will usually follow the custom of identifying $P$ with the poset. If $x$ and $y$ are elements of a poset $P$, $x<y$, and there is no
element $z\in P$ such that $x<z<y$, then we say that $x$ is covered by $y$ (or $y$ covers $x$), and we write $x\prec y$ (or $y\succ x$) \cite{Crawley73}.

A lattice is a poset $P$ such that any two of its elements have the greatest lower bound or meet denoted by $x\wedge y$, and the least upper bound or join denoted by $x\vee y$. A lattice $L$ is complete when each of its subsets $X$ has the greatest lower bound and the least upper bound in $L$. A bounded lattice $L$ is a lattice that possesses the bottom element $0$ and the top element $1$. Note that each complete lattice is a bounded lattice with $\wedge \emptyset=\vee L=1$ and $\vee \emptyset=\wedge L=0$ \cite{Birkhoff73}. Every lattice $L$ in this article with $\mid L\mid\geq 2$.

An element that covers the least element $0$ of a lattice $L$ will be referred to as an atom of $L$. An atomistic lattice is a lattice $L$ in which every element is a join of atoms, and hence of the atoms which it contains. We say that an element $q$ in a lattice $L$ is join-irreducible (resp. meet-irreducible) if, for all $x, y\in L$, $q=x\vee y$ (resp. $q=x\wedge y$) implies $q=x$ or $q=y$. Furthermore, we say that an element $q$ in a lattice is bi-irreducible if $q$ is both join-irreducible and meet-irreducible. The set of join-irreducible elements, the set of atoms of $L$ and the set of bi-irreducible elements will be denoted by $J(L)$, $A(L)$ and $B(L)$, respectively \cite{Birkhoff73}. We say that an element $q$ in a complete lattice $L$ is completely join-irreducible if, for 
every subset $S$ of $L$, $q=\vee S$ implies $q\in S$. Certainly every completely join-irreducible element is join-irreducible \cite{Crawley73}.

We say that a lattice $L$ is modular if, for all $a, b, c\in L$, $$a\geq b \mbox{ implies } a\wedge (b\vee c)=b\vee (a\wedge c).$$
Alternatively, a lattice is modular if and only if it satisfies the following identity:
$$(x\vee y)\wedge ((x\wedge y)\vee z)=(x\wedge y)\vee ((x\vee y)\wedge z).$$
Observe that modularity is a self-dual lattice property\cite{Crawley73}.
\begin{theorem}[\cite{Crawley73, Birkhoff73}]\label{The1}
A lattice $L$ is modular if and only if $L$ has no non-modular five-element sublattice $N_{5}$.
\end{theorem}

In the following, we shall recall some definitions and facts about t-norms and pseudo-t-norms \cite{Wangxue,Wang02, Ma, De}.
\begin{definition}[\cite{Wangxue,Wang02}]\label{De1}
\emph{A binary operation $T: L^2\rightarrow L$ with $L$ a bounded lattice is called a pseudo-t-norm if, for all $a, b, c\in L$, it satisfies \\
(T1) $T(1, a)=a$ and $T(0, a)=0$;\\
(T2) $b\leq c$ implies $T(a, b)\leq T(a, c)$. }
\end{definition}

In what follows, we always assume that the pseudo-t-norm $T$ is commutative, i.e.,  $T(a, b)=T(b, a)$ for any $a, b\in L$. Thus obviously $T(a,1)=a$, $T(a,0)=0$ and $T(a, b)\leq a\wedge b$.

\begin{definition}[\cite{De}]\label{De2}
\emph{A binary operation $T: L^2\rightarrow L$ with $L$ a bounded lattice is called a t-norm if, for all $a, b, c\in L$, it satisfies \\
(a) $T(1, a)=a$;\\
(b) if $b\leq c$ then $T(a, b)\leq T(a, c)$;\\
(c) $T(a, b)=T(b, a)$; \\
(d) $T(T(a, b), c)=T(a, T(b, c))$. }
\end{definition}

Let $L$ be a bounded lattice. A pseudo-t-norm (resp. t-norm) $T: L^2\rightarrow L$ is said to be $\vee$-distributive if, for all $a, b, c\in L$,
\begin{equation*}
 T(a, b\vee c)=T(a,b)\vee T(a, c).
\end{equation*}
A pseudo-t-norm (resp. t-norm) $T: L^2\rightarrow L$ is said to be $\wedge$-distributive if, for all $a, b, c\in L$,
\begin{equation*}
T(a, b\wedge c)=T(a,b)\wedge T(a, c).
\end{equation*}

Let $L$ be a complete lattice. A pseudo-t-norm (resp. t-norm) $T: L^2\rightarrow L$ is said to be left-continuous if, for all $a\in L, \emptyset \neq S\subseteq L$,
\begin{equation}\label{Eq1}
 T(a, \vee S)=\bigvee_{s\in S} T(a,s).
\end{equation}
A pseudo-t-norm (resp. t-norm) $T: L^2\rightarrow L$ is said to be right-continuous if, for all $a\in L, \emptyset \neq S\subseteq L$,
\begin{equation}\label{Eq2}
T(a, \wedge S)=\bigwedge_{s\in S} T(a,s).
\end{equation}
A pseudo-t-norm (resp. t-norm) $T: L^2\rightarrow L$ is continuous if and only if $T$ satisfies both Eqs. \eqref{Eq1} and \eqref{Eq2}.
\begin{remark}\label{remark3.1}
	\emph{For every finite lattice, the left-continuity of a binary operation $T: L^2\rightarrow L$ is equivalent to the $\vee$-distributivity, obviously.}
\end{remark}
\section{Some properties of left-continuous t-norms}
In this section, we give a necessary and sufficient condition for the existence of a left-continuous t-norm on the ordinal sum of two disjoint complete lattices.

In \cite{Birkhoff73}, Birkhoff provided a definition for building the ordinal sum $L_1\oplus L_2$ of two disjoint lattices $L_1$ and $L_2$.
\begin{definition}[\cite{Birkhoff73}]\label{Define1}
\emph{Let $L_1$ and $L_2$ be any two disjoint lattices. The ordinal sum $L_1\oplus L_2$ of $L_1$ and $L_2$ is
the set of all $a\in L_1$ and $b\in L_2$; $a< b$ for all $a\in L_1$ and $b\in L_2$; the relations $a\leq a_1$ and $b\leq b_1$ ($a, a_1\in L_1$; $b, b_1\in L_2$) have unchanged meanings.}
\end{definition}

For complete lattices, the diagram of $L_1\oplus L_2$ can be constructed by laying the diagram of $L_2$ above the diagram of $L_1$ and drawing a line from the least element $0_{L_2}$ of $L_2$ to the greatest element $1_{L_1}$ of $L_1$. In \cite{Gratz}, Gr\"{a}tzer and Wehrung provides a definition for building the glued sum $L_1\dotplus L_2$ of two disjoint complete lattices $L_1$ and $L_2$.

\begin{definition}[\cite{Gratz}]\label{Define2}
\emph{For two disjoint complete lattices $L_1$ and $L_2$. The glued sum $L_1\dotplus L_2$ of $L_1$ and $L_2$ is
defined by putting $L_2$ atop $L_1$ and identifying $1_{L_1}$ with $0_{L_2}$.}
\end{definition}

It is easily seen from Definitions \ref{Define1} and \ref{Define2} that both the ordinal sum and the glued sum of complete lattices $L_1$ and $L_2$ are complete lattices, respectively, 
and these two operations are not commutative in general, i.e., $L_1\oplus L_2\neq L_2\oplus L_1$ and $L_1\dotplus L_2\neq L_2\dotplus L_1$.

Drawing upon Definitions \ref{Define1} and \ref{Define2}, the subsequent observation can be made.

\begin{remark}\label{REMARK}
\emph{Let $L_1$ and $L_2$ be two disjoint complete lattices and $L=L_1\oplus L_2$. Set $L_{\ast}=[0_{L_1}, 0_{L_2}]_{L}$. Then $L=L_{\ast}\dotplus L_2$, 
$1_{L_{\ast}}=0_{L_2}$ and $1_{L_{\ast}}\in J(L_{\ast})\subseteq J(L)$. One can verify that the weakest t-norm 
$$T_{D}(x,y)=\begin{cases}
x\wedge y, &\mbox{if } 1_{L_{\ast}}\in\{x, y\},\\
0, & \mbox{otherwise}.
\end{cases}$$
on $L_{\ast}$ is left-continuous since $1_{L_{\ast}}$ is a completely join-irreducible element in $L_{\ast}$.}
\end{remark}

\begin{proposition}\label{propo1}
Let $L_1, L_2$ be two disjoint complete lattices. If there exist two left-continuous t-norms $T_1$ and $T_2$ on $L_1$ and $L_2$, respectively, then there exists a left-continuous t-norm on $L_1\dotplus L_2$.
\end{proposition}
\proof
Define a binary operation $T: (L_1\dotplus L_2)^2\rightarrow L_1\dotplus L_2$ by
\begin{equation}\label{equation01}
T(x,y)=
\begin{cases}
T_1(x,y), &\mbox{if } x, y\in L_1,\\
T_2(x,y), &\mbox{if } x, y\in L_2,\\
x\wedge y,& \mbox{otherwise}.
\end{cases}
\end{equation}
Proposition 5.2 in \cite{Sam} tells us that $T$ is a t-norm on the complete lattice $L_1\dotplus L_2$.
Now, we shall verify that $T$ is left-continuous by distinguishing three cases as follows.

Case 1. Let $\emptyset\neq S\subseteq L_1$. If $x\in L_1$ then $$T(x, \vee S)=T_1(x, \vee S)=\bigvee_{s\in S}T_1(x, s)=\bigvee_{s\in S}T(x, s)$$ since $T_1$ is a left-continuous t-norm on $L_1$ and $\vee S\in L_1$. If $x\in L_2$ then $$T(x, \vee S)=x\wedge \vee S=\vee S=\bigvee_{s\in S} x\wedge s=\bigvee_{s\in S}T(x, s).$$

Case 2. Let $\emptyset\neq S\subseteq L_2$. If $x\in L_1$ then $$T(x, \vee S)=x\wedge \vee S=x=\bigvee_{s\in S}x\wedge s=\bigvee_{s\in S}T(x, s)$$ since $\vee S\in L_2$ and $x\leq s$ for any $s\in S$. 
If $x\in L_2$ then
$$T(x, \vee S)=T_2(x, \vee S)=\bigvee_{s\in S}T_2(x, s)=\bigvee_{s\in S}T(x, s)$$ since $T_2$ is a left-continuous t-norm on $L_2$ and $\vee S\in L_2$.

Case 3. Let $S=S_1\cup S_2$ with $\emptyset\neq S_1\subseteq L_1$ and $\emptyset\neq S_2\subseteq L_2$. If $x\in L_1$ then $$T(x, \vee S)=T(x, \vee S_2)=x\wedge \vee S_2=x=\bigvee_{s\in S_2} x\wedge s=\bigvee_{s\in S_2} T(x, s)$$ since 
$\vee S=\vee S_2$ and $x\leq s$ for any $s\in S_2$. Note that $\bigvee_{s\in S_1} T(x, s)=\bigvee_{s\in S_1} T_1(x, s)\leq x$.
Thus $$\bigvee_{s\in S} T(x, s)=\bigvee_{s\in S_1} T(x, s) \vee \bigvee_{s\in S_2} T(x, s)=\bigvee_{s\in S_1} T_1(x, s) \vee \bigvee_{s\in S_2} T(x, s)=\bigvee_{s\in S_2} T(x, s),$$ which means that $T(x, \vee S)=\bigvee_{s\in S} T(x, s)$. If $x\in L_2$ then $$T(x, \vee S)=T(x, \vee S_2)=T_2(x, \vee S_2)=\bigvee_{s\in S_2} T_2(x, s)$$ since $T_2$ is a left-continuous t-norm on $L_2$ and $\vee S=\vee S_2$. From Case 1, we have that $T(x, \vee S_1)=\bigvee_{s\in S_{1}}T(x,s)=\vee S_1$. Note that $\vee S_1\leq \bigvee_{s\in S_2} T_2(x, s)$. Hence,
$$T(x, \vee S)=\vee S_1\vee \bigvee_{s\in S_2} T_2(x, s)=\bigvee_{s\in S_1} T(x, s) \vee \bigvee_{s\in S_2} T(x, s)=\bigvee_{s\in S} T(x, s).$$
\endproof

\begin{remark}\label{remark01}
\emph{Let $T_1$ and $T_2$ be two t-noms on two disjoint complete lattices $L_1$ and $L_2$, respectively. If $T$ is defined by Eq. \eqref{equation01}, then $T$ is a left-continuous t-norm on $L_1\dotplus L_2$ if and only if $T_1$ and $T_2$ are left-continuous t-noms on $L_1$ and $L_2$, respectively.}
\end{remark}

\begin{theorem}\label{TTT1}
Let $L_1, L_2$ be two disjoint complete lattices. Then there exists a left-continuous t-norm on $L_1\oplus L_2$ if and only if
there exists a left-continuous t-norm on $L_2$.
\end{theorem}
\proof
From Proposition \ref{propo1} and Remark \ref{REMARK}, the sufficiency is obvious. Now, suppose that $T$ is a left-continuous t-norm on $L_1\oplus L_2$. 
Then, we define a binary operation $T_{\ast}: L_2^2 \rightarrow L_2$ by
\begin{equation}\label{EEE}
T_{\ast}(x,y)=
\begin{cases}
0_{L_2}, &\mbox{if } T(x, y)\in L_1,\\
T(x,y), &\mbox{if } T(x, y)\in L_2.
\end{cases}
\end{equation} 
Evidently, based on Eq. \eqref{EEE} and the fact that $T$ is a t-norm on $L_1\oplus L_2$, $T_{\ast}$ is commutative and has the neutral element $1_{L_2}$.

Let $x\in L_2$ and $S$ be a non-empty set of $L_2$. If $T(x, \vee S)\in L_1$, then we can conclude that $T(x, y)\in L_1$ for all $y\in S$ since $y\leq \vee S$ for every $y\in S$ and $T$ is a monotonically increasing function. Thus, by Eq. \eqref{EEE}, we have 
\begin{equation}\label{EEE1}
T_{\ast}(x,\vee S)=0_{L_2}=\bigvee_{y\in S}T_{\ast}(x,y).  
\end{equation}
 
Now, assume that $T(x, \vee S)\in L_2$. Define $\mathcal{A}=\{y\in S| T(x, y)\in L_1\}$. We assert that $\mathcal{A}\neq S$. 
Suppose, for the sake of contradiction, that $\mathcal{A}=S$. Then, $$T(x, \vee S)=\bigvee_{y\in S}T(x, y)\leq 1_{L_1}<0_{L_2},$$
which contradicts the given condition that $T(x, \vee S)\in L_2$. 

If $\mathcal{A}=\emptyset$, then 
\begin{equation}\label{EEE2}
T_{\ast}(x,\vee S)=T(x, \vee S)=\bigvee_{y\in S}T(x,y)=\bigvee_{y\in S}T_{\ast}(x,y).  
\end{equation} 

If $\mathcal{A}\neq\emptyset$ then $\emptyset \neq \mathcal{A} \subsetneq S$ is obvious. 
Note that $$\bigvee_{y\in \mathcal{A}}T(x,y)\leq 1_{L_1}< 0_{L_2}\leq \bigvee_{y\in S\setminus \mathcal{A}}T(x,y).$$
So, by Eq. \eqref{EEE} and $T(x, \vee S)\in L_2$,
\begin{align*}
T_{\ast}(x,\vee S)&=T(x,\vee S)\\
&=\bigvee_{y\in \mathcal{A}}T(x,y) \vee \bigvee_{y\in S\setminus \mathcal{A}}T(x,y)\\
&=\bigvee_{y\in S\setminus \mathcal{A}}T(x,y)\\
&=0_{L_2}\vee \bigvee_{y\in S\setminus \mathcal{A}}T(x,y)\\
&=0_{L_2}\vee \bigvee_{y\in S\setminus \mathcal{A}}T_{\ast}(x,y)\\
&=\bigvee_{y\in \mathcal{A}}T_{\ast}(x,y) \vee \bigvee_{y\in S\setminus \mathcal{A}}T_{\ast}(x,y).
\end{align*}
This implies that 
\begin{equation}\label{EEE3}
T_{\ast}(x,\vee S)=\bigvee_{y\in S}T_{\ast}(x,y).  
\end{equation} 

Consequently, by Eqs. \eqref{EEE1}, \eqref{EEE2} and \eqref{EEE3}, $T_{\ast}$ is an infinite $\vee$-distributive operation. 

Next, we merely need to prove that $T_{\ast}$ satisfies the associative law. 

Since $T(0_{L_2},y)\leq 0_{L_2}$ for all $y\in L_2$, from Eq. \eqref{EEE}, it immediately follows that
\begin{equation}\label{EEE4}T_{\ast}(0_{L_2}, y)=0_{L_2}\end{equation} 
for each $y\in L_2$. Let $x, y, z\in L_2$. Then we will analyze the associative law with respect to two distinct cases.

Case 1. Suppose $T_{\ast}(x, T_{\ast}(y, z))\neq 0_{L_2}$. By Eq. \eqref{EEE4}, it follows that both $T_{\ast}(y, z)\neq 0_{L_2}$  and $x\neq 0_{L_2}$.
Applying Eq. \eqref{EEE}, we obtain $T_{\ast}(x, T_{\ast}(y, z))=T(x, T(y, z))$. Since $T$ satisfies the associative law, we have $T(x, T(y, z))=T(T(x, y),z)$. This 
implies that $T_{\ast}(x, T_{\ast}(y, z))=T(T(x, y),z)>0_{L_2}$. Thus, we can infer that $T(x, y)> 0_{L_2}$. Consequently, applying Eq. \eqref{EEE} again, we 
conclude that  \begin{equation*}\label{EEE5}T_{\ast}(x, T_{\ast}(y, z))=T(T(x, y),z)=T_{\ast}(T_{\ast}(x, y),z),\end{equation*} which verifies the associative law for $T_{\ast}$ in this case.

Case 2. Assume that $T_{\ast}(x, T_{\ast}(y, z))=0_{L_2}$. Similar to the reasoning in Case 1, we can show that if $T_{\ast}(T_{\ast}(x, y),z)\neq 0_{L_2}$,
then according to the definition of $T_{\ast}$ and the properties of $T$, we get $T_{\ast}(x, T_{\ast}(y, z))\neq 0_{L_2}$, a contradiction. 
Therefore, $$T_{\ast}(x, T_{\ast}(y, z))=0_{L_2}=T_{\ast}(T_{\ast}(x, y),z),$$ verifying the associativity of $T_{\ast}$ for this case.

Combining Cases 1 and 2, we conclude that $T_{\ast}$ satisfies the associative law.

\endproof

\section{Structure of $1$-distributive lattices}
In this section we characterize a complete atomistic lattice on which there is a left-continuous pseudo-t-norm by terms of $1$-distributive lattices. We also describe the forbidden structures of a finite modular lattice that is a $1$-distributive lattice.   
\begin{definition}\label{Defi1}
\emph{Let $L$ be a lattice with the top element $1$. An element $c$ of $L$ is called $1$-distributive if,
for any $a,b\in L$, $a\vee b=1$ implies $c\wedge (a\vee b)=(c\wedge a)\vee (c\wedge b)$. Moreover, a lattice $L$ is called $1$-distributive if every element of $L$ is $1$-distributive.}
\end{definition}

Obviously, every bounded distributive lattice is $1$-distributive. Nevertheless, the converse is not true generally.
In addition, for any bounded lattice $L$, if $1\in J(L)$ then $L$ is $1$-distributive.
Obviously, both $S_{7,2}$ and $S_{7,2}^{*}$ shown in Figure 1 are two $1$-distributive lattices that are non-distributive, in which 
the top element $1$ is not a join-irreducible element.

\par\noindent\vskip50pt
\begin{minipage}{11pc}
\setlength{\unitlength}{0.75pt}\begin{picture}(600,240)
\put(140,80){\circle{6}}\put(136,68){\makebox(0,0)[l]
{\footnotesize $0$}}
\put(140,160){\circle{6}}\put(125,168){\makebox(0,0)[l]
{\footnotesize $b$}}
\put(140,200){\circle{6}}\put(126,200){\makebox(0,0)[l]
{\footnotesize $u$}}
\put(180,160){\circle{6}}\put(185,158){\makebox(0,0)[l]
{\footnotesize $c$}}
\put(180,240){\circle{6}}\put(177,252){\makebox(0,0)[l]
{\footnotesize $1$}}
\put(220,200){\circle{6}}\put(225,198){\makebox(0,0)[l]
{\footnotesize $v$}}
\put(100,160){\circle{6}}\put(85,158){\makebox(0,0)[l]
{\footnotesize $d$}}
\put(100,120){\circle{6}}\put(85,118){\makebox(0,0)[l]
{\footnotesize $m$}}
\put(180,120){\circle{6}}\put(185,118){\makebox(0,0)[l]
{\footnotesize $n$}}
\put(140,163){\line(0,1){34}}
\put(100,123){\line(0,1){34}}
\put(102,122){\line(1,1){36}}
\put(178,122){\line(-1,1){36}}
\put(218,202){\line(-1,1){36}}
\put(102,162){\line(1,1){36}}
\put(142,202){\line(1,1){36}}
\put(178,162){\line(-1,1){36}}
\put(138,82){\line(-1,1){36}}
\put(142,82){\line(1,1){36}}
\put(182,162){\line(1,1){36}}
\put(180,123){\line(0,1){35}}
\put(130,45){$S_{7,2}$}

\put(340,80){\circle{6}}\put(336,68){\makebox(0,0)[l]
{\footnotesize $0$}}
\put(340,120){\circle{6}}\put(325,118){\makebox(0,0)[l]
{\footnotesize $b$}}
\put(340,200){\circle{6}}\put(326,200){\makebox(0,0)[l]
{\footnotesize $u$}}
\put(380,160){\circle{6}}\put(385,158){\makebox(0,0)[l]
{\footnotesize $c$}}
\put(380,240){\circle{6}}\put(377,252){\makebox(0,0)[l]
{\footnotesize $1$}}
\put(420,200){\circle{6}}\put(425,198){\makebox(0,0)[l]
{\footnotesize $v$}}
\put(300,160){\circle{6}}\put(285,158){\makebox(0,0)[l]
{\footnotesize $d$}}
\put(300,120){\circle{6}}\put(285,118){\makebox(0,0)[l]
{\footnotesize $m$}}
\put(380,120){\circle{6}}\put(385,118){\makebox(0,0)[l]
{\footnotesize $n$}}
\put(340,83){\line(0,1){34}}
\put(300,123){\line(0,1){34}}
\put(342,122){\line(1,1){36}}
\put(338,122){\line(-1,1){36}}
\put(418,202){\line(-1,1){36}}
\put(302,162){\line(1,1){36}}
\put(342,202){\line(1,1){36}}
\put(378,162){\line(-1,1){36}}
\put(338,82){\line(-1,1){36}}
\put(342,82){\line(1,1){36}}
\put(382,162){\line(1,1){36}}
\put(380,123){\line(0,1){35}}
\put(333,45){$S_{7,2}^{*}$}
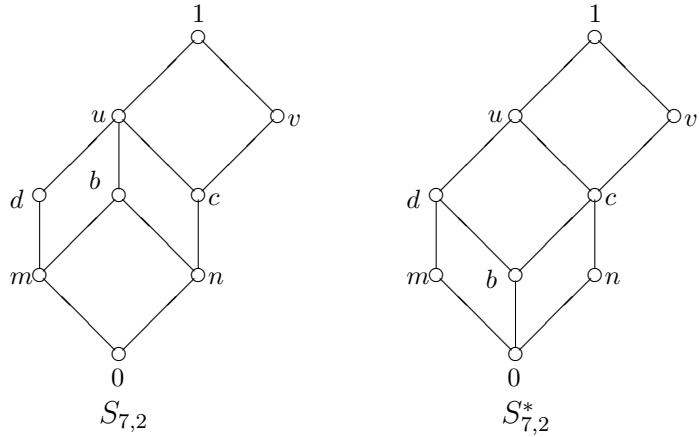
\put(160,20){ Figure 1 Two $1$-distributive lattices.}
\end{picture}
\end{minipage}

\begin{proposition}\label{pr001}
Let $L_1$ and $L_2$ be two disjoint complete lattices. Then $L_1 \dotplus L_2$ is $1$-distributive if and only if 
$L_2$ is $1$-distributive.
\end{proposition}
\proof
Let $L_1 \dotplus L_2$ be $1$-distributive. It is straightforward to see that  $L_2$ is $1$-distributive.
Now, assume that $L_2$ is $1$-distributive. Let $a, b, c\in L_1 \dotplus L_2$ with $a\vee b=1_{L_1 \dotplus L_2}$. 
Evidently, $a$ and $b$ cannot both belong to $L_1$. Thus there are two cases as follows.

Case 1. Let $a, b\in L_2$. If $c\in L_2$ then $c\wedge (a\vee b)=(c\wedge a)\vee (c\wedge b)$ since $L_2$ is $1$-distributive. 
If $c\in L_1$ then $c\leq a$ and $c\leq b$, which yields that $c\wedge (a\vee b)=(c\wedge a)\vee (c\wedge b)$.

Case 2. Let one of $a$ and $b$ is in $L_1$ and the other is in $L_2$. Without loss of generality, let $a\in L_1$ and $b\in L_2$. 
Then $a\vee b=b=1_{L_1 \dotplus L_2}$  and $c\wedge a\leq c\wedge b$ since $a\leq b$. Therefore, $c\wedge (a\vee b)=c\wedge b=(c\wedge a)\vee (c\wedge b)$.

In conclusion, $L_1 \dotplus L_2$ is $1$-distributive.
\endproof

\begin{remark}\label{re001}
\emph{Based on Remark \ref{REMARK} and Proposition \ref{pr001}, it can be deduced that for any two disjoint complete lattices $L_1$ and $L_2$, $L_1 \oplus L_2$ is $1$-distributive if and only if 
$L_2$ is $1$-distributive.}
\end{remark}

\begin{proposition}\label{Prop1}
Let $L$ be a complete atomistic lattice. Then $L$ is $1$-distributive if and only if $L$ is a Boolean lattice.
\end{proposition}
\proof
If $L$ is a Boolean lattice then it is $1$-distributive obviously. 

Conversely, if $L$ is not a Boolean lattice
then $L$ is not isomorphic to the power set of $A(L)$. As a result, there exist two elements $a, b\in L$ such that $a\vee b=1$
while $A(a)\cup A(b)\neq A(L)$. Suppose that $c\in A(L)\setminus (A(a)\cup A(b))$. Then
$c\wedge (a\vee b)=c\neq 0=(c\wedge a)\vee (c\wedge b)$, which means that $L$ is not a $1$-distributive lattice.
\endproof

\begin{remark}\label{rre}
\emph{Proposition \ref{Prop1} indicates that if $L$ is a complete atomistic lattice, then a $1$-distributive lattice $L$ is distributive.}
\end{remark}

\begin{proposition}\label{pr002}
Let $L$ be a bounded lattice. If there exists a $\vee$-distributive pseudo-t-norm on $L$ then $L$ is $1$-distributive.
\end{proposition}
\proof
Let $T: L^2\rightarrow L$ be a $\vee$-distributive pseudo-t-norm. Then $T(c, a\vee b)=T(c,a)\vee T(c,b)$ for any $a, b,c\in L$.
Assume that $a\vee b=1$. Then $T(c, a\vee b)=T(c,1)=c=T(c, a)\vee T(c, b)\leq (c\wedge a)\vee (c\wedge b)\leq c\wedge (a\vee b)=c$.
Consequently, $c\wedge (a\vee b)=(c\wedge a)\vee (c\wedge b)$ for all $c\in L$ whenever $a\vee b=1$. This shows that $L$ is $1$-distributive.
\endproof

\begin{remark}\label{rre1}
\emph{Since every t-norm is a pseudo-t-norm and the left-continuity property implies $\vee$-distributivity, if there exists either a 
left-continuous pseudo-t-norm or a left-continuous t-norm on lattice $L$ then by Proposition \ref{pr002}, $L$ is $1$-distributive.}
\end{remark}

\begin{remark}\label{rre0001}
\emph{Generally, the inverse of Proposition \ref{pr002} does not hold.}
\end{remark}
\begin{example}\label{example1}
\emph{Consider the bounded $1$-distributive lattice $S_{7,2}$ depicted in Figure 1. Let $T$ be a $\vee$-distributive pseudo-t-norm defined on $S_{7,2}$. 
Then $$c=T(c,1)=T(c, d\vee v)=T(c,d)\vee T(c, v)=0\vee T(c,v)=T(c, v)$$ since $T(c, d)\leq c\wedge d=0$. 
Moreover, by the monotonicity of $T$, $c=T(c,v)\leq T(u,v)\leq u\wedge v=c$, i.e., $T(u,v)=c$. 
On the other hand, based on $T(n,v)\leq n\wedge v=n$ and $T(d,v)\leq d\wedge v=0$, we obtain that $$c=T(u, v)=T(d\vee n,v)=T(d, v)\vee T(n, v)=T(n, v)\leq n< c,$$ a contradiction.
Consequently, for a bounded $1$-distributive lattice $L$, there is no a $\vee$-distributive pseudo-t-norm on $S_{7,2}$.}
\end{example}

\begin{theorem}\label{Theo1}
Let $L$ be a complete atomistic lattice with the top element $1$. Then the following are equivalent.\\
\emph{(a)} $L$ is $1$-distributive.\\
\emph{(b)} $L\cong 2^{A(L)}$.\\
\emph{(c)} There exists a continuous t-norm on $L$.\\
\emph{(d)} There exists a left-continuous t-norm on $L$.\\
\emph{(e)} There exists a left-continuous pseudo-t-norm on $L$.\\
\end{theorem}
\proof
By Proposition \ref{Prop1}, (a) $\Rightarrow$  (b). Let $T_{M}$ be such that $T_{M}(x, y)=x\wedge y$ for any $x, y\in L$. Then (b) $\Rightarrow$ (c) since $T_M$ is a continuous t-norm on a
Boolean lattice. Clearly, (c)
$\Rightarrow$ (d) and (d) $\Rightarrow$ (e). According to Remark \ref{rre1}, (e) $\Rightarrow$ (a).
\endproof

\begin{definition}\label{Defi2}
\emph{A sublattice of a lattice $L$ with the top element $1$ is called $1$-sublattice, if it contains the top element $1$.}
\end{definition}

\begin{theorem}\label{Theo2}
A finite modular lattice $L$ with the top element $1$ is $1$-distributive if and only if it does not contain $1$-sublattices $M_3, M_{3,2}$ and $M_{3,4}$
whose Hasse diagrams are presented in Figure 2.

\par\noindent\vskip50pt
\begin{minipage}{11pc}
\setlength{\unitlength}{0.75pt}\begin{picture}(600,200)
\put(60,120){\circle{6}}
\put(58,132){\makebox(0,0)[l]{\footnotesize$1$}}
\put(60,40){\circle{6}}
\put(20,80){\circle{6}}
\put(60,80){\circle{6}}
\put(100,80){\circle{6}}
\put(60,43){\line(0,1){34}}
\put(60,83){\line(0,1){34}}
\put(58,42){\line(-1,1){36}}
\put(62,42){\line(1,1){36}}
\put(22,82){\line(1,1){36}}
\put(98,82){\line(-1,1){36}}
\put(50,15){$M_3$}

\put(260,120){\circle{6}}
\put(260,40){\circle{6}}
\put(220,80){\circle{6}}
\put(180,120){\circle{6}}
\put(220,160){\circle{6}}
\put(218,172){\makebox(0,0)[l]{\footnotesize$1$}}
\put(260,80){\circle{6}}
\put(300,80){\circle{6}}
\put(260,43){\line(0,1){34}}
\put(260,83){\line(0,1){34}}
\put(258,42){\line(-1,1){36}}
\put(218,82){\line(-1,1){36}}
\put(262,42){\line(1,1){36}}
\put(182,122){\line(1,1){36}}
\put(222,82){\line(1,1){36}}
\put(298,82){\line(-1,1){36}}
\put(258,122){\line(-1,1){36}}
\put(250,15){$M_{3,2}$}

\put(460,120){\circle{6}}
\put(460,40){\circle{6}}
\put(420,80){\circle{6}}
\put(380,120){\circle{6}}
\put(420,160){\circle{6}}
\put(458,212){\makebox(0,0)[l]{\footnotesize$1$}}
\put(460,80){\circle{6}}
\put(500,80){\circle{6}}
\put(540,120){\circle{6}}
\put(500,160){\circle{6}}
\put(460,200){\circle{6}}

\put(460,43){\line(0,1){34}}
\put(460,83){\line(0,1){34}}
\put(458,42){\line(-1,1){36}}
\put(418,82){\line(-1,1){36}}
\put(462,42){\line(1,1){36}}
\put(502,82){\line(1,1){36}}
\put(382,122){\line(1,1){36}}
\put(422,162){\line(1,1){36}}
\put(422,82){\line(1,1){36}}
\put(462,122){\line(1,1){36}}
\put(498,82){\line(-1,1){36}}
\put(538,122){\line(-1,1){36}}
\put(498,162){\line(-1,1){36}}
\put(458,122){\line(-1,1){36}}
\put(450,15){$M_{3,4}$}
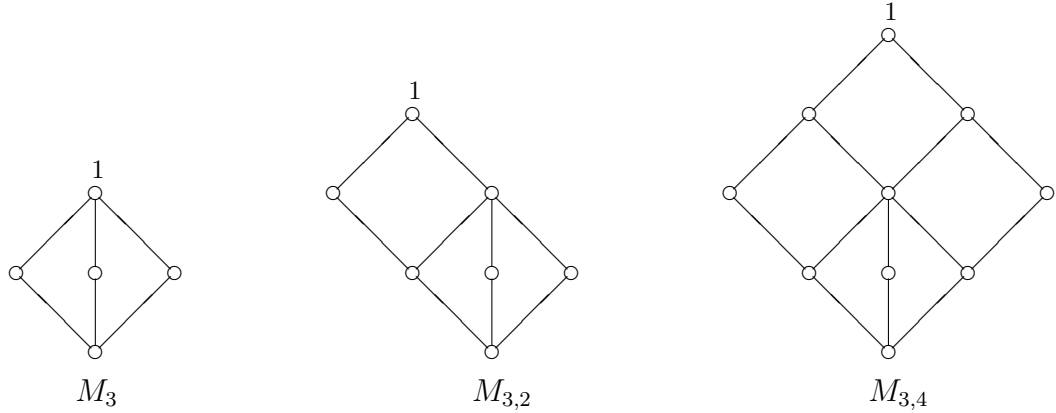
\put(110,-5){\emph{Figure 2 Three non-$1$-distributive modular lattices}}
\end{picture}
\end{minipage}
\end{theorem}
\proof
The necessity is clear. Now, let $L$ denote a finite non-$1$-distributive modular lattice. 
Set $$\mathcal{S}=\{x\in L|x \mbox{ is a } \mbox{non-}1\mbox{-distributive element}\}.$$
Obviously, $\mathcal{S}\neq \emptyset$ by Definition \ref{Defi1}. Since $L$ is finite, there
exists a maximal element in $\mathcal{S}$. Let $c$ be a maximal element in $\mathcal{S}$ and 
$$\mathcal{R}=\{(x,y)\in L^{2}| x\vee y=1, c>(c\wedge x)\vee (c\wedge y)\}.$$
Evidently, $\mathcal{R}\neq \emptyset$. Define $\leq$ on $\mathcal{R}$ as follows: 

for any $(x,y), (u,v)\in \mathcal{R}$, $(x,y)\leq(u,v)$ if and only if $x\leq u$ and $y\leq v$.\\
It is easy to see that $(\mathcal{R}, \leq)$ is a finite poset. Let $(a,b)$ be a maximal element of $\mathcal{R}$. Then we have the following two statements.\\
($\diamond$) If $a\wedge c=b\wedge c$ then $a\wedge b=a\wedge c=b\wedge c$.

In reality, if $a\wedge b\neq a\wedge c$ and $a\wedge b\neq b\wedge c$, then $a\wedge b> a\wedge b\wedge c =a\wedge c=b\wedge c$, which implies that $a\wedge b\parallel  c$. In fact, if $a\wedge b\geq c$, then $(a\wedge c)\vee (b\wedge c)=c$, contrary to $(a,b)\in \mathcal{R}$. If $a\wedge b\leq c$, then $a\wedge b\leq a\wedge c$, contrary to the fact that $a\wedge b> a\wedge c$.
Thus $a\wedge b\parallel  c$, which means that $c< (a\wedge b)\vee c$ and $a\wedge b< (a\wedge b)\vee c$. Furthermore, by modularity we know that
\begin{equation*}\label{Equa2}
\{[(a\wedge b)\vee c]\wedge a\}\vee \{[(a\wedge b)\vee c]\wedge b\} =[(a\wedge b)\vee (a\wedge c)]\vee [(a\wedge b)\vee (b\wedge c)]=a\wedge b
\end{equation*}
since $a\wedge b>a\wedge c=b\wedge c$, i.e., $\{[(a\wedge b)\vee c]\wedge a\}\vee \{[(a\wedge b)\vee c]\wedge b\} < (a\wedge b)\vee c$. This follows that $(a\wedge b)\vee c \in \mathcal{S}$, contrary to the fact that $c$ is a maximal element in $\mathcal{S}$.\\
($\triangle$) $a\parallel c$ and $b\parallel c$.

If $a\geq c$ then $(c\wedge a) \vee (c\wedge b)=c\vee (c\wedge b)=c$. This contradicts the fact that $(a, b)\in \mathcal{R}$. If $a\leq c$, then we have 
$$(c\wedge a)\vee (c\wedge b)=a\vee (c\wedge b)=(a\vee c)\wedge (a\vee b)=c\wedge (a\vee b)=c$$
since $L$ is modular and $a\vee b=1$, contrary to $(a,b)\in \mathcal{R}$. Therefore, $a\parallel c$. Analogously, it can be proved that $b\parallel c$.

The rest of the proof is structured by considering three distinct cases.

(I) $a\vee c=b\vee c$. 

In this case, we assert that $a\wedge c=b\wedge c$. In fact, if $a\wedge c< b\wedge c$ then $a\wedge c=a\wedge b\wedge c<b\wedge c$. This implies that $a\ngeqslant b\wedge c$, and thus $a< a\vee (b\wedge c)$. 
Then $(a\vee(b\wedge c), b)\geq (a, b)$ and $(a\vee(b\wedge c), b)\neq (a, b)$, hence $(a\vee(b\wedge c), b)\notin \mathcal{R}$ since $(a,b)$ is a maximal element of $\mathcal{R}$.
Therefore, by modularity and ($\triangle$), we have
\begin{equation}\label{Equa1}
c=[(a\vee (b\wedge c))\wedge c]\vee (b\wedge c)=[(b\wedge c)\vee (a\wedge c)]\vee (b\wedge c)=b\wedge c<c,
\end{equation}
a contradiction. Consequently, $a\wedge c\nless b\wedge c$. It can be similarly proved that $b\wedge c\nless a\wedge c$.
Now, suppose that $a\wedge c\parallel b\wedge c$. If $a=a\vee (b\wedge c)$ then $a\geq b\wedge c$. Thus $a\wedge c\geq b\wedge c$, a contradiction.
Hence $a<a\vee (b\wedge c)$. Analogous to the above discussion, the last inequality will lead to Eq. \eqref{Equa1}, a contradiction.
As a result, we conclude that $a\wedge c=b\wedge c$, which together with ($\diamond$) means that $a\wedge b=a\wedge c=b\wedge c$. Consequently, there exists 
a $1$-sublattice $M_{3}=\{a, b, c, a\wedge b\wedge c, 1\}$ in $L$.

(II) Either $a\vee c>b\vee c$ or $a\vee c<b\vee c$. 

In the following, we just prove this theorem when $a\vee c>b\vee c$. One can analogously verify the theorem when $a\vee c<b\vee c$.  Evidently, in this case we have $a\vee c=a\vee b\vee c=1$. We claim that $a\wedge c=b\wedge c$. Otherwise, there are three cases as follows.

Case 1. If $a\wedge c< b\wedge c$ then $a\ngeq b\wedge c$. This follows that $a\parallel b\wedge c$. Indeed, if $a< b\wedge c$, then 
$a\vee c=c \leq b\vee c$, a contradiction. Thus $a\parallel b\wedge c$. Hence $a<a\vee (b\wedge c)$. Analogous to the above discussion, the last inequality  will lead to Eq. \eqref{Equa1}, a contradiction.

Case 2. If $a\wedge c>b\wedge c$, then similar to Case 1 this will lead to a contradiction.

Case 3. If $a\wedge c\parallel b\wedge c$ then $a\vee (b\wedge c)> a$. In fact, if $a\vee (b\wedge c)=a$ then
$b\wedge c\leq a$. This implies that $b\wedge c\leq a\wedge c$, a contradiction. Thus similar to the proof of Case 1 the inequality $a\vee (b\wedge c)> a$ will lead to a contradiction.

Cases 1, 2 and 3 mean that $a\wedge c=b\wedge c$. Thus by ($\diamond$), $a\wedge b=a\wedge c=b\wedge c$. 
On the one hand, we claim that $a\wedge (b\vee c)\neq a\wedge b\wedge c$. In fact, if $a\wedge (b\vee c)= a\wedge b\wedge c$ then the
five elements $a\wedge b\wedge c, c, b\vee c, a$ and $1$ form an $N_5$-isomorphic sublattice, contradicting the modularity of $L$.
Thus $a\wedge b\wedge c< a\wedge (b\vee c)$. On the other hand, we assert that $a\wedge (b\vee c)< a$. Indeed, if $a\wedge (b\vee c)=a$ then $a\leq b\vee c$. 
This yields that $a\vee c\leq b\vee c$, a contradiction. Therefore, $a\wedge b\wedge c< a\wedge (b\vee c)<a$. 
Meanwhile, $[a\wedge (b\vee c)]\wedge c=a\wedge c=a\wedge b\wedge c$ and
$[a\wedge (b\vee c)]\wedge b=a\wedge b=a\wedge b\wedge c$.
Moreover, by modularity and $b\vee c< a\vee c=a\vee b=1$, we know that $[a\wedge (b\vee c)]\vee c=(b\vee c)\wedge (a\vee c)=b\vee c$ and
$[a\wedge (b\vee c)]\vee b=(b\vee c)\wedge (a\vee b)=b\vee c$. Consequently, there exists a $1$-sublattice
$M_{3,2}=\{a, b, c, a\wedge b\wedge c, a\wedge (b\vee c), b\vee c, 1\}$ in $L$.

(III) $a\vee c\parallel b\vee c$. Similar to the proofs of Cases 1, 2 and 3 in (II), we can show $a\wedge c=b\wedge c$.
Combining with ($\diamond$), we get $a\wedge b=a\wedge c= b\wedge c$.
We claim that $b\wedge (a\vee c)\neq b\wedge c$. Otherwise, the elements $a\wedge b\wedge c, b, a, a\vee c$ and $1$ form an $N_5$-isomorphic 
sublattice, a contradiction. Similarly, we can prove that $a\wedge (b\vee c)\neq a\wedge c$.
Moreover, using modularity and $a\vee b=1$, we have that
$$c\vee [b\wedge (a\vee c)]=(a\vee c)\wedge (b\vee c),$$
$$c\vee [a\wedge (b\vee c)]=(b\vee c)\wedge (a\vee c),$$
$$a\vee [b\wedge (a\vee c)]=(a\vee c)\wedge (a\vee b)=a\vee c,$$
$$b\vee [a\wedge (b\vee c)]=(b\vee c)\wedge (a\vee b)=b\vee c$$
and
\begin{align*}
[a\wedge (b\vee c)]\vee [b\wedge (a\vee c)]
&=\{[a\wedge (b\vee c)]\vee (a\vee c)\}\wedge \{[a\wedge (b\vee c)]\vee b\}\\
 &=(a\vee c)\wedge [(b\vee c)\wedge (a\vee b)]\\
 &=(a\vee c)\wedge (b\vee c).
\end{align*}
Consequently, there exists a $1$-sublattice
$M_{3,4}=\{a, b, c, a\wedge b\wedge c, a\wedge (b\vee c), b\wedge (a\vee c), a\vee c, b\vee c, (a\vee c)\wedge (b\vee c), 1\}$ in $L$.
\endproof

Theorem \ref{Theo2} and Proposition \ref{pr002} deduce the following lemma.

\begin{lemma}\label{lemm001}
Let $L$ be a finite modular lattice with the top element $1$. If there exists a $\vee$-distributive pseudo-t-norm (resp. t-norm) on $L$, then $L$ 
does not contain $1$-sublattices $M_3, M_{3,2}$ and $M_{3,4}$.
\end{lemma}

\section{Left-continuous pseudo-t-norm}

In this section we characterize a kind of finite planar modular lattices on which there is a left-continuous pseudo-t-norm by terms of $1$-distributive lattices.

A lattice $L$ is planar if its Hasse diagram 
can be drawn for $L$ in which none of the straight line segments intersect\cite{Kelly}. 

Let $C$ with the top element $1_C$ and the bottom element $0_C$ and $D$ with the top element $1_D$ and the bottom element $0_D$  be two finite chains. A left corner of $L=C\times D$ is defined as follows. 
Note that $(1_C,0_D)\in l$ and it is a bi-irreducible element where $l$ is the left boundary 
chain of $L$. Removing $(1_C,0_D)$, we get the lattice $L_1$. Having defined a planar distributive lattice $L_{n-1}$ with 
a left boundary chain $l_{n-1}$, pick a bi-irreducible element $a\in l_{n-1}$. Define $L_{n}=L_{n-1}-\{a\}$ with a left
boundary chain $l_{n}$. The left corner of $L$ is $L-L_{n}$. We can define right corners similarly. 

\begin{lemma}[\cite{Gratz07}]\label{LmmaKell}
A finite, planar, distributive lattice can be obtained from the direct product of two finite chains by removing a left and a right corner.
\end{lemma}

It is well known that a planar modular lattice $S^{+}$ can be derived from a planar distributive lattice $S$ by adding ``eyes" to covering 
squares-making covering squares into covering $M_{n}$-s, see the lattice $S^{+}$ in Figure 3 for an example; 
the elements ``eyes" in $S^{+}$ are black-filled \cite{Gratz10}.

\par\noindent\vskip50pt
\begin{minipage}{11pc}
\setlength{\unitlength}{0.75pt}\begin{picture}(600,240)
\put(140,120){\circle{6}}
\put(140,200){\circle{6}}
\put(60,200){\circle{6}}
\put(100,240){\circle{6}}
\put(140,40){\circle{6}}
\put(100,80){\circle{6}}
\put(60,120){\circle{6}}
\put(100,160){\circle{6}}
\put(98,252){\makebox(0,0)[l]{\footnotesize$1$}}
\put(138,30){\makebox(0,0)[l]{\footnotesize$0$}}
\put(180,80){\circle{6}}
\put(180,160){\circle{6}}
\put(138,42){\line(-1,1){36}}
\put(98,82){\line(-1,1){36}}
\put(142,42){\line(1,1){36}}
\put(62,122){\line(1,1){36}}
\put(102,82){\line(1,1){36}}
\put(178,82){\line(-1,1){36}}
\put(138,122){\line(-1,1){36}}
\put(142,122){\line(1,1){36}}
\put(98,162){\line(-1,1){36}}
\put(138,202){\line(-1,1){36}}
\put(178,162){\line(-1,1){36}}
\put(102,162){\line(1,1){36}}
\put(62,202){\line(1,1){36}}
\put(130,10){$S$}
\put(400,120){\circle{6}}
\put(400,200){\circle{6}}
\put(320,200){\circle{6}}
\put(360,240){\circle{6}}
\put(400,40){\circle{6}}
\put(360,80){\circle{6}}
\put(380,76){$\bullet$}
\put(412,76){$\bullet$}

\put(356,198){$\bullet$}
\put(320,120){\circle{6}}
\put(360,160){\circle{6}}
\put(358,252){\makebox(0,0)[l]{\footnotesize$1$}}
\put(398,30){\makebox(0,0)[l]{\footnotesize$0$}}
\put(440,80){\circle{6}}
\put(440,160){\circle{6}}
\put(360,200){\line(0,1){36}}
\put(360,163){\line(0,1){36}}
\put(398,42){\line(-1,1){36}}
\put(358,82){\line(-1,1){36}}
\put(402,42){\line(1,1){36}}
\put(322,122){\line(1,1){36}}
\put(362,82){\line(1,1){36}}
\put(438,82){\line(-1,1){36}}
\put(398,122){\line(-1,1){36}}
\put(398,42){\line(-2,5){14}}
\put(402,42){\line(2,5){14}}
\put(402,122){\line(1,1){36}}
\put(358,162){\line(-1,1){36}}
\put(398,202){\line(-1,1){36}}
\put(438,162){\line(-1,1){36}}
\put(362,162){\line(1,1){36}}
\put(322,202){\line(1,1){36}}
\put(382,78){\line(2,5){15.5}}
\put(416,78){\line(-2,5){15.5}}
\put(390,10){$S^{+}$}
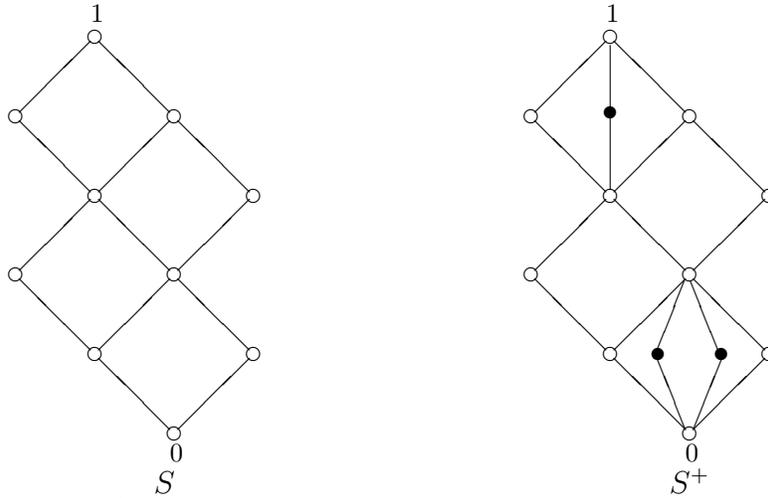
\put(45,-5){Figure 3 A planar distributive and a planar modular lattice}
\end{picture}
\end{minipage}

\begin{definition}\label{de5.1}
\emph{A rectangular lattice is defined as a planar lattice $L$ with exactly two bi-irreducible elements, which are complementary and distinct from $0$ and $1$, on 
the boundary of $L$.} 
\end{definition}

Obviously, a finite rectangular distributive lattice is isomorphic to the lattice $C\times D$, where $C$ and $D$ represent two finite chains. 
Furthermore, every finite rectangular modular lattice can be derived by adding ``eyes" to a finite rectangular distributive lattice.  It can be verified that 
both the lattices $S$ and $S^{+}$ presented in Figure 3 are not rectangular. 

\begin{theorem}\label{theorem5.1}
Let $L$ be a finite rectangular modular lattice with the top element $1$. Then the following are equivalent.\\
\emph{(i)} $L$ is $1$-distributive.\\
\emph{(ii)} $L$ is distributive.\\
\emph{(iii)} There exists a continuous t-norm on $L$.\\
\emph{(iv)} There exists a left-continuous t-norm on $L$.\\
\emph{(v)} There exists a left-continuous pseudo-t-norm on $L$. 
\end{theorem}
\proof
(i) $\Rightarrow$ (ii) by Theorem \ref{Theo2} and Definition \ref{de5.1}.

(ii) $\Rightarrow$ (iii) is straightforward because $T_M$ is a continuous t-norm on $L$.

(iii) $\Rightarrow$ (iv) is trivial.

(iv) $\Rightarrow$ (v) is trivial.

(v) $\Rightarrow$ (i) by Proposition \ref{pr002} and Remark \ref{rre1}.
\endproof

\begin{remark}
\emph{Observe that for a finite planar lattice $L$ with the top element $1$ and the bottom element $0$, if $0\in B(L)$ (resp. $1\in B(L)$) then $L$ is a sum of a finite chain and a finite planar sublattice. 
Thus, in the following, we always assume that $\{0, 1\}\nsubseteq B(L)$.}
\end{remark}

\begin{theorem}\label{theorem5.2}
Let $L$ be a finite planar modular lattice with the top element $1$ and exactly one bi-irreducible element on either the left boundary chain $C_l$ or the right boundary chain $C_r$. 
Then $L$ is $1$-distributive if and only if there exists a $\vee$-distributive pseudo-t-norm on $L$.
\end{theorem}
\proof
The sufficiency is straightforward. Now, assume that $L$ is a $1$-distributive planar modular lattice. According to Definition \ref{de5.1}, if both the chains $C_l$ and $C_r$ have exactly one bi-irreducible element, then $L$ is a rectangular modular lattice. Hence, by 
Theorem \ref{theorem5.1}, there exists a $\vee$-distributive pseudo-t-norm on $L$. Next, we only need to consider the case that only one of the $C_l$ and $C_r$ 
has exactly one bi-irreducible element. Without loss of generality, assume $C_l$ has exactly one such element $a$ while $b_1< \cdots < b_n=b$ are 
bi-irreducible elements on the right boundary $C_r$. By Theorem \ref{Theo2}, $[a\wedge b, 1]$ forms a 
rectangular distributive sublattice of $L$ since $L$ is $1$-distributive, i.e., $[a\wedge b, 1]=[a\wedge b, a] \times [a\wedge b, b]$. 
Define a binary operation $T: L^2\rightarrow L$ as follows. 
\begin{equation}\label{eq5.1}
T(x,y)=\begin{cases}
0, & (x,y)\in [0,b)\times H_1\cup H_1\times[0,b),\\
a\wedge x\wedge y, & (x,y)\in H_2\times H_2,\\
x\wedge y, & \mbox{otherwise},
\end{cases} 
\end{equation}
where $H_1=L-[b,1]$ and $H_2=H_1-[0,b)$. 
In what follows, we prove that $T$ is a $\vee$-distributive pseudo-t-norm on $L$.
Evidently, $T$ is commutative with the neutral element $1$. Because of $\vee$-distributivity implying 
monotonicity, according to Definition \ref{De1}, we only need to prove 
that $T$ satisfies that $T(x,y\vee z)=T(x, y)\vee T(x,z)$ for any $x, y, z\in L$. 
Observe first that Eq. \eqref{eq5.1} tells us that 
\begin{equation}\label{eq5.11}
T(x, y)\leq x\wedge y
\end{equation}
for any $x, y\in L$. We distinguish three cases as below.

I. If $y\vee z\in [0, b)$, then $y\in [0, b)$ and  $z\in [0, b)$. If $x\in H_1$ then $$T(x,y\vee z)=0=0\vee 0=T(x, y)\vee T(x,z).$$ 
If $x\in [b,1]$ then $$T(x, y\vee z)=x\wedge (y\vee z)=y\vee z=(x\wedge y)\vee (x\wedge z)=T(x, y)\vee T(x,z).$$

II. If $y\vee z\in [b, 1]$, then, by the structure of $L$, either $y$ or $z$ belongs to $[b, 1]$. Without loss of generality, assume $y\in [b, 1]$. Now, we 
consider three distinct cases.

Case 1. If $x\in [0, b)$, then $T(x, y)=x\wedge y=x$. Since Eq. \eqref{eq5.11}  implies that $T(x, z)\leq x\wedge z\leq x$, 
we have that $$T(x, y\vee z)=x\wedge (y\vee z)=x=x\vee T(x, z)=T(x, y)\vee T(x,z).$$   

Case 2. If $x\in [b, 1]$, then in the case $z\in [0, b)$, $$T(x, y\vee z)=x\wedge (y\vee z)=x\wedge y=T(x, y)=T(x, y)\vee T(x,z)$$ since 
$T(x,y)=x\wedge y\geq b> z=x\wedge z=T(x, z)$. If $z\in L-[0, b)$ then $x, y, z \in [a\wedge b, 1]$. Since $[a\wedge b, 1]$ forms a 
rectangular distributive sublattice of $L$, we know that  
\begin{equation}\label{equation5.2}
T(x, y\vee z)=x\wedge (y\vee z)=(x\wedge y)\vee (x\wedge z)=T(x, y)\vee T(x,z).
\end{equation}

Case 3. If $x\in H_2$, then in the case $z\in [0, b)$, $$T(x, y\vee z)=T(x, y)=T(x, y)\vee T(x,z)$$ 
since $T(x, z)=0$ and $y> z$. If $z\in [b,1]$, then $$T(x, y\vee z)=T(x, y)\vee T(x,z)$$ as in the proof of Eq. \eqref{equation5.2}. 
In the case $z\in H_2$, we know that $x, y$ and $z$ belong to the distributive sublattice $[a\wedge b, 1]$. Then
\begin{align*}
	x\wedge (y\vee z)&=(x\wedge y)\vee (x\wedge z)\\
	&=[(x\wedge y)\vee (x\wedge z)]\wedge (x\vee a) \ \ \ \mbox{ since } (x\wedge y)\vee (x\wedge z)\leq x\vee a\\
	&=[(x\wedge y)\vee (x\wedge z)]\wedge [(x\vee a)\wedge (a\vee b)] \ \ \ \mbox{ since } a\vee b=1\\
	&=[(x\wedge y)\vee (x\wedge z)]\wedge [(x\wedge b)\vee a]\\
	&\leq [(x\wedge y)\vee (x\wedge z)]\wedge [(x\wedge y)\vee a]  \ \ \ \mbox{ since } y\in [b,1]\\
	&=(x\wedge y)\vee [(x\wedge z)\wedge a]\\
	&\leq(x\wedge y)\vee (x\wedge z)\\
	&=x\wedge (y\vee z),
\end{align*}
which implies that $x\wedge (y\vee z)=(x\wedge y)\vee (x\wedge z\wedge a)$. Thus $$T(x, y\vee z)=x\wedge (y\vee z)=(x\wedge y)\vee (x\wedge z\wedge a)=T(x, y)\vee T(x,z).$$

III. If $y\vee z\in H_2$, then either $y\in H_2$ or $z\in H_2$, say $y\in H_2$. 
By the structure of $L$, we have that  
\begin{equation}\label{eq5.12}
z\in H_1.                   
\end{equation}
Then the other discussion is carried out in the following three cases.

Case a. If $x\in H_2$, then by Eq. \eqref{eq5.12}, in the case $z\in [0, b)$ we have $$T(x, y\vee z)=a\wedge x\wedge (y\vee z)\leq a\wedge x\wedge (y\vee b)
=(a\wedge x\wedge y)\vee (a\wedge x\wedge b)=a\wedge x\wedge y$$ since  $a\wedge x, y, b\in [a\wedge b, 1]$ 
and $[a\wedge b, 1]$ is a distributive sublattice. Thus $T(x, y\vee z)=a\wedge x\wedge y$. This together with  $T(x, z)=0$ implies that 
$$T(x, y\vee z)=a\wedge x\wedge y=T(x, y)=T(x, y)\vee T(x,z).$$ 
In the case $z\in H_2$, $x, y$ and $z$ belong to the distributive sublattice $[a\wedge b, 1]$. Hence 
$$T(x, y\vee z)=a\wedge x\wedge (y\vee z)= (a\wedge x\wedge y)\vee (a\wedge x\wedge z)=T(x, y)\vee T(x,z). $$

Case b. If $x\in [0, b)$, then, by Eq. \eqref{eq5.12}, $$T(x, y\vee z)=0=0\vee 0=T(x, y)\vee T(x,z). $$

Case c. If $x\in [b, 1]$, then by Eq. \eqref{eq5.12}, in the case $z\in H_2$, similar to the proof of Eq. \eqref{equation5.2} we have 
$$T(x, y\vee z)=T(x, y)\vee T(x,z).$$
If $z\in[0, b)$, then $$T(x, y\vee z)=x\wedge (y\vee z)=(x\wedge z)\vee (x\wedge y)=T(x, y)\vee T(x,z)$$ since $L$ is modular and $z< x$. 

In summary, the binary operation $T$ defined by Eq.\eqref{eq5.1} is a $\vee$-distributive pseudo-t-norm on $L$.
\endproof

\begin{example}\label{exam5.1}
\emph{Consider the lattice $L$ as shown in Figure 4.}
\end{example}

\par\noindent\vskip50pt
\begin{minipage}{11pc}
\setlength{\unitlength}{0.75pt}
\begin{picture}(600,200)
\put(240,120){\circle{6}}
\put(240,200){\circle{6}}
\put(240,40){\circle{6}}
\put(200,80){\circle{6}}
\put(188,76){\makebox(0,0)[l]{\footnotesize$a$}}
\put(220,76){$\bullet$}
\put(227,78){\makebox(0,0)[l]{\footnotesize$b$}}
\put(252,76){$\bullet$}
\put(260,78){\makebox(0,0)[l]{\footnotesize$c$}}
\put(160,120){\circle{6}}
\put(148,118){\makebox(0,0)[l]{\footnotesize$f$}}
\put(248,118){\makebox(0,0)[l]{\footnotesize$e$}}
\put(200,160){\circle{6}}
\put(188,160){\makebox(0,0)[l]{\footnotesize$g$}}
\put(238,212){\makebox(0,0)[l]{\footnotesize$1$}}
\put(238,30){\makebox(0,0)[l]{\footnotesize$0$}}
\put(280,80){\circle{6}}
\put(286,78){\makebox(0,0)[l]{\footnotesize$d$}}
\put(280,160){\circle{6}}
\put(287,158){\makebox(0,0)[l]{\footnotesize$h$}}
\put(238,42){\line(-1,1){36}}
\put(198,82){\line(-1,1){36}}
\put(242,42){\line(1,1){36}}
\put(162,122){\line(1,1){36}}
\put(202,82){\line(1,1){36}}
\put(278,82){\line(-1,1){36}}
\put(238,122){\line(-1,1){36}}
\put(238,42){\line(-2,5){14}}
\put(242,42){\line(2,5){14}}
\put(242,122){\line(1,1){36}}
\put(278,162){\line(-1,1){36}}
\put(202,162){\line(1,1){36}}
\put(222,78){\line(2,5){15.5}}
\put(256,78){\line(-2,5){15.5}}
\put(230,10){$L$}
\put(125,-5){Figure 4 A planar modular lattice}
\end{picture}
\end{minipage}\\

It can be verified that $L$ satisfies all the conditions of Theorem \ref{theorem5.2}. Then, applying Eq.\eqref{eq5.1}, we construct the $\vee$-distributive 
pseudo-t-norm $T$ on $L$ as shown by Table \ref{Tab:01}. However, $T$ is not a t-norm since $T(T(f, g), h)=T(f,h)=a\neq 0=T(f,e)=T(f,T(g,h))$.
\begin{table}[htpp]
\centering
\caption{The $\vee$-distributive pseudo-t-norm $T$ on $L$}
\label{Tab:01}
\begin{tabular}{c|cccccccccc}

 $T$ & $0$ & $a$ & $b$ & $c$ & $d$ & $e$ & $f$ & $g$ & $h$ & $1$ \\
 \hline
$0$ & $0$ & $0$ & $0$ & $0$ & $0$ & $0$ & $0$ & $0$ & $0$ & $0$\\
$a$ & $0$ & $0$ & $0$ & $0$ & $0$ & $0$ & $0$ & $0$ & $a$ & $a$\\
$b$ & $0$ & $0$ & $0$ & $0$ & $0$ & $0$ & $0$ & $0$ & $b$ & $b$\\
$c$ & $0$ & $0$ & $0$ & $0$ & $0$ & $0$ & $0$ & $0$ & $c$ & $c$\\
$d$ & $0$ & $0$ & $0$ & $0$ & $0$ & $0$ & $0$ & $0$ & $d$ & $d$\\
$e$ & $0$ & $0$ & $0$ & $0$ & $0$ & $0$ & $0$ & $0$ & $e$ & $e$\\
$f$ & $0$ & $0$ & $0$ & $0$ & $0$ & $0$ & $f$ & $f$ & $a$ & $f$\\
$g$ & $0$ & $0$ & $0$ & $0$ & $0$ & $0$ & $f$ & $f$ & $e$ & $g$\\
$h$ & $0$ & $a$ & $b$ & $c$ & $d$ & $e$ & $a$ & $e$ & $h$ & $h$\\
$1$ & $0$ & $a$ & $b$ & $c$ & $d$ & $e$ & $f$ & $g$ & $h$ & $1$\\
\end{tabular}
\end{table}

\section{Concluding remarks}
In this article, we successfully established some equivalent relations among the continuity of t-norms, 
the $1$-distributivity of lattices and distributive lattices. By using three finite sublattices, 
we obtained a necessary and sufficient condition for a finite modular lattice being $1$-distributive lattices. Meanwhile, we equivalently 
characterized a kind of planar $1$-distributive modular lattices on which there is a left-continuous pseudo-t-norm. Overall, this article contributes 
to a deeper understanding of the connection between the algebraic structures of lattices and the continuity of pseudo-t-norm even if it just partly answer the two questions raised in Section 1.

\end{document}